\input amstex
\input epsf
\documentstyle{amsppt}
% \Monograph
\pagewidth{450pt}
\pageheight{700pt}
\loadbold

\define\R#1{\ifcase#1
\or1%   R1
\or1%   R2
\or1%   R3
\or4%   R4
\or2%   R5
\or3%   R6
\or50%  R7
\or5%   R8
\or6%   R9
\or7%   R10
\or8%   R11
\or9%   R12
\or10%  R13
\or11%  R14
\or15%  R15
\or16%  R16
\or17%  R17
\or18%  R18
\or19%  R19
\or51%  R20
\fi}

\define\Zerner{Jackson }
\define\Sobolev{Kolmogorov }

%% \overfullrule=0pt

%% the parameters above should be properly set before printing

\define\GMI{\operatorname{GMI}}
\define\HCP{\operatorname{HCP}}
\define\LS{\operatorname{\L S}}

\define\JP{\operatorname{JP}}

\define\capa{\operatorname{cap}}

\define\diam{\operatorname{diam}}
\define\dist{\operatorname{dist}}
\define\inte{\operatorname{int}}

\define\xb#1{B_{i_1,\dots,i_{#1}}}
\define\xz#1{z_{i_1,\dots,i_{#1}}}

\define\diagram#1{{\normallineskip=6pt \normalbaselineskip=0pt \matrix #1 \endmatrix}}

\define\harr#1#2{\smash{\mathop{\hbox to .5in{\rightarrowfill}}\limits^{\scriptstyle#1}_{\scriptstyle#2}}}
\define\Harr#1#2{\smash{\mathop{\hbox to .5in{\leftarrowfill}}\limits^{\scriptstyle#1}_{\scriptstyle#2}}}
\define\hHarr#1#2{\smash{\mathop{\hbox to .5in{\leftarrowfill\hskip -2pt\rightarrowfill}}\limits^{\scriptstyle#1}_{\scriptstyle#2}}}
\define\varr#1#2{\llap{$\scriptstyle #1$}\left\downarrow\vcenter to .5in{}\right.\rlap{$\scriptstyle #2$}}
\define\Varr#1#2{\llap{$\scriptstyle #1$}\left\uparrow\vcenter to .5in{}\right.\rlap{$\scriptstyle #2$}}
\define\vVarr#1#2{\llap{$\scriptstyle #1$}\left\updownarrow\vcenter to .5in{}\right.\rlap{$\scriptstyle #2$}}
\define\varrne#1#2{\llap{$\scriptstyle #1$}\nearrow\vcenter to .5in{}\rlap{$\scriptstyle #2$}}
\define\varrnw#1#2{\llap{$\scriptstyle #1$}\nwarrow\vcenter to .5in{}\rlap{$\scriptstyle #2$}}
\define\varrsw#1#2{\llap{$\scriptstyle #1$}\swarrow\vcenter to .5in{}\rlap{$\scriptstyle #2$}}
\define\varrse#1#2{\llap{$\scriptstyle #1$}\searrow\vcenter to .5in{}\rlap{$\scriptstyle #2$}}
\define\varrswne#1#2{\llap{$\scriptstyle #1$}\swarrow\nearrow\vcenter to .5in{}\rlap{$\scriptstyle #2$}}
\define\Diagram#1{{\normallineskip=8pt \normalbaselineskip=0pt \matrix #1 \endmatrix}}

\catcode`\@=11
\def\Eqalign#1{\null\,\vcenter{\openup\jot\m@th
\ialign{\strut\hfil$\displaystyle{##}$&
$\displaystyle{{}##}$\hfil&& \quad $\displaystyle{##}$\hfil&
\qquad $\displaystyle{{}##}$\hfil\crcr #1\crcr}}\,} \catcode`\@=12

\catcode`\@=11
\def\EqaligN#1{\null\,\vcenter{\openup\jot\m@th
\ialign{\strut\hfil$\displaystyle{##}$&
$\displaystyle{{}##}$\hfil&& \, $\displaystyle{##}$\hfil& \qquad
$\displaystyle{{}##}$\hfil\crcr #1\crcr}}\,} \catcode`\@=12

\input cyracc.def
\catcode`\@=11 \font@\tencyr=wncyr10 \font@\eightcyr=wncyr8
\catcode`\@=13 \addto\tenpoint{}
\addto\eightpoint{}

\topmatter

\title
Jackson's inequality in the complex plane and the \L ojasiewicz-Siciak inequality of Green's function
\endtitle

\rightheadtext{Jackson's inequality and Green's function}
%\NoRunningHeads

\author
Leokadia BIALAS-CIEZ and Raimondo EGGINK
\endauthor

\affil
Institute of Mathematics, Jagiellonian University,
% \address
ul. \L ojasiewicza 6, 30-348 Krak\'ow, Poland,
% \endaddress
% \email
e-mail: leokadia.bialas-ciez\@im.uj.edu.pl,
% \endemail
tel. +48 12 664 66 34,
fax +48 12 664 66 74.
\endaffil

\abstract We prove a generalization of Dunham Jackson's famous
approximation inequality to the case of compact sets in the
complex plane admitting both upper and lower bounds for their
Green's functions, i.e.\ the well known H\"older Continuity
Property ($\HCP$) and the less known but crucial \L
ojasiewicz-Siciak inequality ($\LS$). Moreover, we show that
($\LS$) is a necessary condition for our Jackson type inequality.
\endabstract

\endtopmatter

\document
\nologo

% \keywords

%\vskip 5mm

%% \subjclass
%% \endsubjclass

\vskip 15mm

\head 1. Introduction and main result.
\endhead

Dunham Jackson's famous inequality which gives some control over
the rate of approximation by polynomials of a fixed function, was
first proved for the segment [-1,1] in 1911 (see \cite{15} and
also \cite{23, sec.5.1; 9, chap.4 sec.6}). There are numerous
results in the literature concerning various generalizations of
this inequality because of their significant role in approximation
theory and in related domains of research. This explains also why
sets admitting Jackson type inequalities are especially useful.
However, it seems that in the complex case this property  was
investigated up to now only for particular classes of sets.

The direct reason for our study of Jackson's inequality was a
result by Bos and Milman  regarding the equivalence of the local
and global Markov inequalities, a Kolmogorov type inequality and
an extension property for $\Cal C^\infty(K)$ functions (see
\cite{7} or L.P. Bos and P.D. Milman, {\it A Geometric
Interpretation and the Equality of Exponents in Markov and
Gagliardo-Nirenberg (Sobolev) Type Inequalities for Singular
Compact Domains}, preprint). The proof is hard and proceeds only
in the real case making essential use of the Jackson inequality in
$\Bbb R^N$. We were intrigued to obtain a corresponding result for
sets in the complex plane because of the intricate
interconnectedness of multiple distinct global and local
properties: Markov inequalities, Kolmogorov type inequalities,
polynomial approximation, extension operators, geometric
properties and, ultimately, the behavior of the Green's function,
i.e.\ L-regularity, H\"older continuity and the \L
ojasiewicz-Siciak inequality. However, a simple adaptation to the
complex case of the proof given by Bos and Milman is not possible.

In a previous paper \cite{5} we showed that the local Markov
property is equivalent to a Kolmogorov type property for any
compact set $E\subset\Bbb C$. In a subsequent paper
(L.\ Bialas-Ciez and R.\ Eggink, {\it Equivalence of the global and
local Markov inequalities in the complex plane}, preprint) we
prove that the Global Markov Inequality GMI (i.e.\ $\|p'\|_E\le M
(\deg p)^k\|p\|_E$ with $k,M>0$ independent of the polynomial $p$)
is equivalent to an extension property for functions of the class
$s(E)$, which can be rapidly approximated by holomorphic
polynomials:
$$s(E) \ := \ \{f\in \Cal C(E)\quad:\quad\forall \, \ell\in \Bbb N \quad
\lim_{n\to\infty}n^\ell\,\dist_E(f,\Cal P_n)=0\},$$ where $\Bbb
N=\{1,2\ldots\}$, $\dist_E(f,\Cal P_n):=\inf\{\|f-p\|_E\ :\
p\in\Cal P_n\}$ is the error of approximating the function $f$ on
the set $E$ by polynomials of degree $n$ or less and $\|\cdot
\|_E$ is the supremum norm on $E$. The latter extension property
requires the existence of an extension, which is bounded together
with its derivatives by the following Jackson norms of the
extended function:
$$|f|_\ell \ := \ \|f\|_E+\sup_{n\in \Bbb N}n^\ell\, \dist_E(f,\Cal P_n)
\ \ \ \ {\text{for}} \ \ \ \ell\in \Bbb N_0=\{0,1,2\ldots\}.$$
Sometimes we will use $|\cdot |_\ell$ also for $\ell\in \Bbb R$,
$\ell\ge 0$.

In the real case contemplated by Bos and Milman, this extension
property implies a Kolmogorov type inequality owing to the fact
that the Jackson norms can easily be estimated by quotient norms.
This follows from the classical Jackson inequality and therefore
we investigated the possibility to generalize this result to the
case of compact sets in the complex plane.

Clearly, a lot of work has been done on various
"J(ackson)-properties" for Jordan arcs, domains and other
continua, where order of approximation is linked to the regularity
of a given function and/or the regularity of the continuum,
see for example \cite{20; 24; 17; 11; 2; 16; 1; 12; 10} and many other
authors referenced therein. However, our research of the
literature leads us to believe that this is not at all the case
for compact sets in general, which may even be totally
disconnected.

One can envisage different possible generalizations of the Jackson inequality,
so we have taken an approach that seems to be best suited to determine a class of sets
for which the global and local Markov inequalities are equivalent.
This allows us to work only with functions that are holomorphic in open neighborhoods
of our compact set and with regular supremum norms in those neighborhoods,
while maintaining optimal control over the constants.

For a compact set $E\subset \Bbb C$, let $\Cal
H^\infty(E):=\left\{\,f\in\Cal C^\infty(\Bbb C)\ :\ \frac{\partial
f}{\partial\bar z}\equiv0\ \text{in some open neighborhood of}\
E\,\right\}$ and $E_\delta:=\left\{z\in\Bbb C\ :\ \dist
(z,E)\le\delta\right\}$. By Taylor's theorem and Cauchy's integral
formula, we can prove that for a closed disc $B\subset \Bbb C$ and
for an arbitrary function $f\in \Cal H^\infty (B_\delta)$ with
some $\delta\in(0,1]$, we have $f_{|B}\in s(B)$ and $$\forall \
\ell \in\Bbb N \quad \quad |f_{|B}|_\ell\ \le \ \left(c \ell
\right)^{\ell+1} \|f\|_{B,\ell +1},$$ where $c$ depends only on
the diameter of $B$. Consequently,
$$\forall \ \ell\ge1\quad  \quad |f_{|B}|_\ell \ \le \
\left(\frac{c\ell}{\delta}\right)^{\ell+1} \|f\|_{B_\delta}.$$ The
Jackson Property defined below is a generalization of the last
inequality.

\definition{Definition 1.1}
A \ compact \ set \ $E\subset\Bbb C$ \ admits \ the \
$\underline{\text{Jackson \ Property}}$ \ $\JP(s)$, \ where
$s\ge1$, \ if \ $\Cal H^\infty(E)_{|E}\subset s(E)$ and
there exist constants $c,v\ge1$ such that
$$
% \exists \, c_0, c_1\ge0\quad
\forall\,\ell\in\Bbb N \quad
\forall\,\delta\in (0,1] \quad \forall \,f\in\Cal H^\infty(E_\delta)\quad:\quad
|f_{|E}|_\ell \ \le \ \left(\frac{c\ell^v}{\delta^s}\right)^{\ell+c} ||f||_{E_\delta}. \tag 1$$
\enddefinition

Note that every closed disc admits $\JP(1)$. Note also that if
$\Cal H^\infty(E)_{|E}\subset s(E)$ then the set $E$ must
obviously be polynomially convex, i.e.\ $E=\hat E$ where $\hat
E:=\{z\in\Bbb C\, :\, \forall n\in \Bbb N \ \forall p\in \Cal P_n
\ |p(z)|\le \|p\|_E\}$ is the polynomial hull of $E$.

\vskip 2mm

The interesting thing is that the Jackson Property defined above
turns out to be intimately connected with the rate of growth of
the Green's function $g_E$ (with logarithmic pole at infinity) of
the unbounded complement of the compact set $E$.

%the \L ojasiewicz-Siciak inequality ($\LS$), which gives a lower
%bound for the Green function of the complement of a given compact
%set

\definition{Definition 1.2}
The set $E$ admits the $\underline{\text{\L ojasiewicz-Siciak
inequality}}$ $\LS(s)$, where $s\ge1$, if
$$\exists \, M>0\quad
\forall\, z\in E_1\quad:\quad g_E(z) \ \ge \ M \dist (z,E)^s.$$
% Note that we extend the Green's function to the entire complex plane
% by putting $g_E(z):=0$ for $z\in\hat E$.
We will write that the set $E$ admits ${\LS}$ if it admits $\LS(s)$
for some $s\ge1$.
\enddefinition

As far as we know, the term \L ojasiewicz-Siciak inequality was
first coined by Gendre, who used it to obtain advanced
approximation results \cite{14} (see also \cite{22}). The
interested reader is referred to \cite{6} for basic information.

We set out (without proofs) the following examples:

$\bullet$ \ if $E$ is a compact set in  $\Bbb{R}$ then $E$ admits
$\LS(1)$,

$\bullet$ \ the set $E := \{z \in \Bbb C : |z -1| \le 1$ or $ |z
+1| \le 1\}$ does not admit $\LS(s)$ for any $s$, \pagebreak

$\bullet$ \ if $E$ is the starlike set $E = E(n) := \left\{z=\, r
\, \exp{\frac {2\pi i  j}n}\in\Bbb C\ :\ 0\le r\le1,\
j=1,\dots,n\right\}$ then $E$ admits $\LS(\frac{n}2)$ whenever
$n\in\Bbb N\setminus\{1\}$,

$\bullet$ \ a simply connected compact set $E\subset\Bbb C$ with
nonempty interior, admits \L S$(s)$ with some $s\ge 1$\linebreak
if and only if its complement to the Riemann sphere is a H\"older
domain, i.e.\ a conformal map \linebreak $\varphi:\{z\in \Bbb C \,
: \, |z|<1\} \rightarrow\hat\Bbb C\setminus E$ \ such that
$\varphi(0)=\infty$ is H\"older continuous in $\{z\in \Bbb C \, :
\, \frac12\le |z|\le 1\}$ with exponent $1/s$.

\vskip 2mm

The {\L}ojasiewicz-Siciak inequality is the opposite of the well
known H\"older Continuity Property ($\HCP$), which gives an upper
bound of the Green's function (see e.g. \cite{8; 3; 19}).

\definition{Definition 1.3}
A compact set $E\!\subset\!\Bbb C$ admits the
$\underline{\text{H\"older Continuity Property}}$ $\HCP(k)$, where
$k\ge1$, if
$$\exists \, M\ge1\quad
\forall \, z\in E_1\quad :\quad g_E(z)\ \le \ M\dist (z,E)^{1/k}.$$
We will write that the set $E$ admits ${\HCP}$ if it admits
$\HCP(k)$ for some $k\ge1$.
\enddefinition

The connection between the Jackson property and the rate of growth
of the Green's function is evidenced by our main result:

\proclaim{Theorem 1.4} Let $s'>s\ge 1$. Any polynomially convex
compact set $E\subset\Bbb C$ admitting $\LS(s)$ and $\HCP$, admits
$\JP(s)$. Moreover, any compact set $E\subset\Bbb C$ admitting
$\JP(s)$, admits $\LS(s')$.
\endproclaim

This finding allowed us to construct an example of a compact set
in the complex plane which admits the Global Markov Inequality,
while it does not admit any Local Markov Property, nor the \L ojasiewicz-Siciak inequality
(L.\ Bialas-Ciez and R.\ Eggink, {\it Equivalence of the global and
local Markov inequalities in the complex plane}, preprint).

This paper is organized as follows. In Section 2 we introduce the
notations used throughout the paper, some of which are standard
while others are more specific to our work. Section 3 contains the
proof of the main result. In Section 4 we give some remarks and
additional results  concerning the Jackson Property. We wrap up
with some open problems.

\vskip 5mm

\head
2. Preliminaries and notations.
\endhead

In our further deliberations we make active use of Siciak's
$\underline{\text{extremal function}}$ for a compact set
$E\subset\Bbb C$ (see \cite{21})
$$\Phi_E(z) \, := \ \limsup_{n\rightarrow\infty}\root n\of{\Phi_n(z)}\quad
\text{for }z\in\Bbb C,$$
where for $n\in\Bbb N$
$$\Phi_n(z) \, = \ \Phi_n(E,z) \, := \ \sup\{|p(z)|\ :\ p\in \Cal P_n\ ,\ \|p\|_E\le 1\}$$
denotes the $n$-th extremal function. It is well known that
$\Phi_E=e^{g_E}$, where $g_E$ stands for the Green's function of
$\Bbb C \setminus \hat E$ with logarithmic pole at infinity. For
convenience we extend $g_E$ to the entire complex plane
by putting $g_E(z):=0$ for all
$z\in\hat E$.

The set $E$ is called $\underline{\text{L-regular}}$ if its
extremal function $\Phi_E$ is continuous on the entire complex
plane. Similarly we speak of regularity in a boundary point
$z_0\in E$ when $\Phi_E$ is continuous at this point. Note that
whenever the cardinality of the set $E$ is bigger than $n$, then
the $n$-th extremal function $\Phi_n$ of $E$ is necessarily
continuous on the entire complex plane. Since the extremal
function $\Phi_E$ is always continuous on $\Bbb C\setminus\hat E$,
\ L-regularity is really determined by the behavior of $\Phi_E$ at
the outer boundary of the set $E$.

Note that both properties $\HCP$ and $\LS$ can be defined
equivalently in terms of Siciak's extremal function instead of
Green's function, because for arbitrary $t>0$ we have
$$1+g_E(z) \ \le \
e^{g_E(z)} \ = \ \Phi_E(z) \ \le \ 1+\frac{e^t-1}t\ g_E(z)$$ for
all $z\in\Bbb C\setminus\hat E$ such that $0\le g_E(z)\le t$.

For a compact set $E\subset\Bbb C$ and $\rho\ge1$ we denote the
level set of the extremal function
$$C(E,\rho) \ := \ \{z\in\Bbb C\ :\ \Phi_E(z)=\rho\} \ \ \{z\in\Bbb C\ :\ g_E(z)=\,\log\rho\}.$$
In order to control the
behavior of the extremal function $\Phi_E$ near the boundary of
$E$ we introduce
$$\aligned
\phi_n(t) \
:= \ \inf_{z\in d E_t}\Phi_n(z),\\
\endaligned$$ for
$n\in\Bbb N_0$ and $t\in[0,\infty)$.
Here and further we denote by
d$E_t$ the set $\{z\in\Bbb C\ :\ \dist(z,E)=t\}$, which may be a
slightly bigger set than just the boundary $\partial E_t$ of
$E_t$. Note that for $t>0$ the function $\phi_n$ is continuous or
equal to $+\infty$. Furthermore, $\phi_n(t)$ is an increasing
function with respect to $n$ and moreover, the maximum principle
for subharmonic functions, applied to the function $\log\Phi_{n}$,
implies that $ \phi_n(t)$ is increasing also with respect to
$t>0$.

For $\delta>0$ we denote by $K(E,\delta)$ a compact neighborhood
constructed as follows. First we cut up the entire complex plane
into closed squares of size $\delta\times\delta$, starting at the
origin of the plane. Next we select all squares having a non-empty
intersection with the set $E$ and by $K(E,\delta)$ we denote the
sum of those squares. Clearly we have $E\subset K(E,\delta)\subset
E_{\delta\sqrt 2}$. Also it is easy to see that the set
$K(E,\delta)$ consists of at most $\left(\frac{\diam
E}\delta+2\right)^2$ squares and therefore the length of its
border $\partial K(E,\delta)$ is definitively less than
$$4\delta\,\left(\frac{\diam E}\delta+2\right)^2 \ = \ \frac{4(\diam E_\delta)^2}\delta.$$

For a compact set $E\subset\Bbb C$ we denote the family of smooth
functions that are $\bar{\partial}$-flat on $E$:
$$\Cal A^\infty(E) \ := \ \left\{\ f\in\Cal C^\infty(\Bbb C)\ :\
\text{the function}\ \frac{\partial f}{\partial\bar z}\ \text{is
flat on}\ E\ \right\},$$ where a function $g\in\Cal C^\infty(\Bbb
C)$ is said to be flat in the point $z_0$ if $D^\alpha g(z_0)=0$
for all $\alpha=(\alpha_1,\alpha_2)\in\Bbb N_0^2$,
$D^\alpha=\frac{\partial^{|\alpha|}}{\partial
z^{\alpha_1}\cdot\partial\bar z^{\alpha_2}}$ and
$|\alpha|=\alpha_1+\alpha_2$. This definition is slightly
different than in \cite{22}, where $\Cal A^\infty(E)$ stood for
functions defined on $E$ only, which will be denoted here as $\Cal
A^\infty(E)_{|E}:=\{f_{|E}\ :\ f\in\Cal A^\infty(E)\}$.

\vskip 5mm

\head 3. Proof of the main result.
\endhead

Our goal in this section is to establish the main result of the
paper, which is a general version of Jackson's inequality in the
complex plane. For a fixed compact set $E\subset\Bbb C$ and
$\zeta\notin E$ we put $f_\zeta(z):=\frac1{\zeta-z}$ for $z$ in
some open neighborhood of $E$ and extend it to a function of class
$\Cal C^\infty(\Bbb C)$ so that $f_\zeta\in\Cal H^\infty(E)$.

\proclaim{Lemma 3.1} For all $\zeta\notin E\subset\Bbb C$ and
$n\in\Bbb N_0$ we have
$$\frac1{\bigl(\dist(\zeta,E)+\diam E\bigr)\ \Phi_{n+1}(\zeta)} \ \le
\ \dist_E(f_\zeta,\Cal P_n) \ \le \ \frac1{\dist(\zeta,E) \
\Phi_{n+1}(\zeta)}.$$
\endproclaim

\demo{Proof} Fix $n\in\Bbb N$ and take an arbitrary polynomial
$q\in\Cal P_{n+1}$ such that $\|q\|_E=1$ and $q(\zeta)\ne 0$.
Define $p(z):=\frac{q(\zeta)-q(z)}{(\zeta-z)\,q(\zeta)}$ so that
$p\in\Cal P_n$. We obtain
$$\dist_E(f_\zeta,\Cal P_n)\le\|f_\zeta-p\|_E\sup_{z\in E}\left|\frac{q(z)}{(\zeta-z)\, q(\zeta)}\right|\le
\frac{\|q\|_E}{|q(\zeta)|\cdot\inf_{z\in E}|\zeta-z|} \frac1{\dist(\zeta,E)\,|q(\zeta)|}.$$
We take the infimum over all
$q\in\Cal P_{n+1}$ to arrive at $\dist_E(f_\zeta,\Cal
P_n)\le\frac1{\dist(\zeta,E)\ \Phi_{n+1}(\zeta)}$.

On the other hand for fixed $n\in\Bbb N$ find $p\in\Cal P_n$ such
that $\dist_E(f_\zeta,\Cal P_n)=\|f_\zeta-p\|_E$. Define
$q(z):=1-(\zeta-z)\,p(z)$ \ so that $q\in\Cal P_{n+1}$. We see
that
$$\|q\|_E= \sup_{z\in E}|1-(\zeta-z)\,p(z)|\le
\sup_{z\in E}|\zeta-z|\cdot\sup_{z\in E}|f_\zeta(z)-p(z)|\le
\bigl(\dist(\zeta,E)+\diam E\bigr) \dist_E(f_\zeta,\Cal P_n)$$ and
hence
$$\Phi_{n+1}(\zeta) \ \ge \
\frac{|q(\zeta)|}{\|q\|_E} \ \ge \
\frac1{\bigl(\dist(\zeta,E)+\diam E\bigr) \dist_E(f_\zeta,\Cal
P_n)}.\qed$$
\enddemo

The next results were inspired by the proof of Runge's theorem
(see e.g. \cite{13, chap.II\S3, chap.III\S1}).

\proclaim{Proposition 3.2} For any compact set $E\subset\Bbb C$,
$0<\delta\le 1$ and $f\in\Cal H^\infty(E_\delta)$ we have
$$\forall \ \frac12\le b<1\quad\forall \ n\in\Bbb N\quad: \quad \dist_E(f,\Cal P_n)\ \le
\
\frac{c\,\|f\|_{E_\delta}}{(1-b)\delta^2\,\phi_{n+1}(b\delta)},$$
where the constant $c:=\frac{28}\pi\left(2+\diam E\right)^2$
depends only on the set $E$.
\endproclaim \pagebreak

\demo{Proof} Fix $\frac12\le b<1$ and $n\in\Bbb N$. If
$\phi_{n+1}(b \delta)=+\infty$ then the set $E$ consists of $n+1$
or less points and $\dist_E(f,\Cal P_n)=0$, which finishes the
proof. Otherwise, find a positive $\widetilde\delta$ such that
$\widetilde\delta\le \frac{(1-b) \delta}{4 \phi_{n+1}(b \delta)}$
and $\frac{(1-b) \delta}{4 \widetilde\delta}$ is an integer. Let
$\Gamma$ be the boundary $\partial K\!\!\left(E_{b
\delta},\frac{1-b}4 \delta\right)$, with proper orientation, and
cut it up into equal intervals $\Gamma_j$, each of length
$\widetilde\delta$, so that $\Gamma=\bigcup_{j}\Gamma_j$, with $j$
running over a finite index set. As $K\!\!\left(E_{b
\delta},\frac{1-b}4 \delta\right)\subset E_{b \delta+\frac{1-b}4
\delta \sqrt2} \subset E_{\frac{1+b}2 \delta}$, we see that
$\Gamma\subset E_{\frac{1+b}2 \delta}\setminus \inte E_{b
\delta}$, \ while for the length of $\Gamma$, denoted $m(\Gamma)$,
we have
$$\sum_j\,\widetilde\delta \ = \ m(\Gamma) \ \le \
\frac{4 (\diam E_\delta)^2}{\frac{1-b}4 \delta} \ \le \ \frac{4\pi
c}{7 (1-b) \delta}.\tag 2$$

For a fixed $z\in E$ and $f\in\Cal H^\infty(E_\delta)$ put
$g_z(\zeta):=\frac{f(\zeta)}{\zeta-z}$, which is a holomorphic
function in an open neighborhood of the set
$E_\delta\setminus\{z\}$. Let $\zeta_0,\zeta_1\in\Gamma_j$ for
some $j$. Then the entire interval $I:=[\zeta_0,\zeta_1]$ lies in
$\Gamma_j$ and of course $\dist(z,I)\ge b \delta$. By Cauchy's
integral formula, for $\zeta\in I$ we have
$$|g_z'(\zeta)| = \left|\frac{f'(\zeta)}{\zeta-z}-
\frac{f(\zeta)}{(\zeta-z)^2}\right|
\le\frac{\|f'\|_{E_{\frac{1+b}2 \delta}}}{b \delta}+
\frac{\|f\|_{E_{\frac{1+b}2 \delta}}}{\left(b \delta\right)^2}\le
\frac{\|f\|_{E_\delta}}{\frac{1-b}2  b \delta^2}+
\frac{\|f\|_{E_\delta}}{\left(b \delta\right)^2} = \frac{(1+b)
\|f\|_{E_\delta}}{(1-b)  b^2 \delta^2}\le \frac{6
\|f\|_{E_\delta}}{(1-b) \delta^2}.$$ This leads us to
$$\left|\frac{f(\zeta_1)}{\zeta_1-z}-\frac{f(\zeta_0)}{\zeta_0-z}\right|
=|g_z(\zeta_1)-g_z(\zeta_0)| \le \int_I |g_z'(\zeta)|\, |d\zeta|
\le\frac{6 \|f\|_{E_\delta}}{(1-b) \delta^2} |\zeta_1-\zeta_0|\le
\frac{6 \|f\|_{E_\delta}}{(1-b) \delta^2} \, \widetilde\delta\le
\frac{6 \|f\|_{E_\delta}}{4\delta \phi_{n+1}(b \delta)}.$$ We now
see that for all $z\in E$, all $j$ and arbitrarily selected points
$\zeta_j\in\Gamma_j$ we have $$ \left| \int_{\Gamma
_j}\frac{f(\zeta )}{\zeta -z}\, d\zeta -
 \int_ {\Gamma_ j}\frac{f(\zeta _j)}{\zeta _j-z}\,
d\zeta\right|\le \int_{\Gamma_j}\left|\frac{f(\zeta)}{\zeta-z}-
\frac{f(\zeta_j)}{\zeta_j-z}\right|\, |d\zeta| \le
\int_{\Gamma_j}\frac{6 \|f\|_{E_\delta}}{4\delta \ \phi_{n+1}(b
\delta)}\, |d\zeta| = \frac{3 \widetilde\delta\ \|f\|_{E_\delta}}
{2 \delta \ \phi_{n+1}(b \delta)}.$$ By summing over $j$ we obtain
$$|f(z)-R(z)| \ = \ \left|\frac1{2\pi i} \int_{\Gamma}\frac{f(\zeta)}{\zeta-z}\, d\zeta-R(z)\right|\
\le \ \sum_j\frac{3 \widetilde\delta \ \|f\|_{E_\delta}} {4\pi
\delta \ \phi_{n+1}(b \delta)} \ = \ \frac{3  m(\Gamma)\
\|f\|_{E_\delta}} {4\pi \delta \ \phi_{n+1}(b \delta)},$$ where we
denote
$$R(z)\ := \ \frac1{2\pi i}\sum_j \int_{\Gamma_j}\frac{f(\zeta_j)}{\zeta_j-z}\, d\zeta
\ = \ \sum_j\frac{c_j}{\zeta_j-z} \ = \ \sum_j c_jf_{\zeta_j}(z),
\ \ \ \ \ \ \ \ c_j:=\frac1{2\pi i} f(\zeta_j)
\int_{\Gamma_j}d\zeta .$$ By the above, we can see that the
rational function $R$ approximates uniformly $f$  on the set $E$
and
$$\|f-R\|_E \ \le \ \frac{3  m(\Gamma)\  \|f\|_{E_\delta}}
{4\pi \delta \ \phi_{n+1}(b \delta)}.\tag 3$$ Simultaneously, by
virtue of Lemma 3.1 and by the minimum principle, we have
$$\dist_E(R,\Cal P_n)\!\le\!\!\sum_j|c_j|  \dist_E(f_{\zeta_j},\Cal P_n)\!\le
\!\!\sum_j\frac{\widetilde\delta\ \|f\|_{E_\delta} }{2\pi
\dist(\zeta_j,E) \Phi_{n+1}(\zeta_j) }
\!\le\!\!\sum_j\frac{\widetilde\delta\ \|f\|_{E_\delta} }{2\pi b
\delta \phi_{n+1}(b \delta)}\! = \!\frac{m(\Gamma)\
\|f\|_{E_\delta} }{2\pi b \delta\phi_{n+1}(b \delta) },$$ because
$\dist(\zeta_j,E)\ge b \delta$. Consequently, from $(2)$ and
$(3)$, since $\frac 1b\le2$, we conclude that
$$\dist_E(f,\Cal P_n) \ \le \
\|f-R\|_E+\dist_E(R,\Cal P_n) \ \le \ \frac{7 m(\Gamma)\
\|f\|_{E_\delta}} {4\pi \delta \ \phi_{n+1}(b \delta)} \ \le \
\frac{c\, \|f\|_{E_\delta}}{(1-b) \delta^2 \, \phi_{n+1}(b
\delta)}.\qed$$
\enddemo

\proclaim{Lemma 3.3} For any L-regular compact set $E\subset\Bbb
C$, $\zeta\in E_1\setminus\hat E$, $1<\rho\le \Phi_E(\zeta)$ and
$n\in\Bbb N_0$ we have
$$\dist_E(f_\zeta,\Cal P_n) \ \le \
\frac{(n+1)\,c}{\dist\bigl(C(E,\rho),E\bigr) \, \dist(\zeta,E)}
\left(\frac{\rho}{\Phi_E(\zeta)}\right)^{n+1},$$ where $c\ge1$
depends only on the set $E$.
\endproclaim \pagebreak

\demo{Proof}
We put
$$d \ := \ \max_{z\in C\left(E,\|\Phi_E\|_{E_1}\right)}\dist(z,E) \ \ge \ 1$$
and $c:=2d+\diam E$. Fix $\zeta\in E_1\setminus\hat E$, $1<\rho\le
\Phi_E(\zeta)$, $n\in\Bbb N_0$ and consider any $\eta\in
E_d\setminus \hat E$. For the Lagrange interpolation polynomial
$L_nf_\eta$ with knots in $n+1$ Fekete extremal points
$\left\{z^{(n)}_j\right\}_{j=0,\dots,n}\subset E$ and
$\omega_n(z):=\prod_{j=0}^n\left(z-z^{(n)}_j\right)$, we have
$$L_nf_\eta(z) \ = \ \frac{\omega_n(\eta)-\omega_n(z)}{\omega_n(\eta) (\eta-z)}.$$
Consequently, applying the properties of the Fekete extremal
points, we see that for all $z\in E$ we have
$$\left|\frac{\omega_n(z)}{\omega_n(\eta)}\right|\!=\!|1\!-\!(\eta-z)
L_nf_\eta(z)|\!\le\! 1+\bigl(d+\diam E\bigr) (n+1)
||f_\eta||_E\!\le 1+ \frac{(d+\diam E)
(n+1)}{\dist(\eta,E)}\!\le\! \frac{(n+1)\,c}{\dist(\eta,E)}.$$

Now put $h_n:=\log|\omega_n|-(n+1)g_E$, which is a harmonic
function on $\Bbb C\setminus\hat E$, bounded in $\hat \Bbb C$. If
 $\eta\in C(E,\rho)\subset E_d\setminus\hat E$ and $z\in E$, we
have
$$h_n(\eta)\ge \log \frac{\dist\left(\eta,E\right)  |\omega_n(z)|}{(n+1)\,c}-(n+1)\log\rho\ge
\log\left(\frac{\dist\bigl(C(E,\rho),E\bigr)
|\omega_n(z)|}{(n+1)\,c}\right)-(n+1) \log\rho.$$ The L-regularity
of the set $E$ leads us to the fact that the level set $C(E,\rho)$
is the boundary of the open domain $\Omega:=\{z\in\Bbb C\ :\
\Phi_E(z)>\rho\}$. Therefore, the minimum principle for harmonic
functions implies that the last inequality holds for all
$\eta\in\bar{\Omega}$, in particular, for $\eta=\zeta$. By the
definition of $h_n$ and since $g_E=\log \Phi_E$, we can easily
obtain
$$\left|\frac{\omega_n(z)}{\omega_n(\zeta)}\right| \ \le \
\frac{(n+1)\,c}{\dist\bigl(C(E,\rho),E\bigr)}
\left(\frac{\rho}{\Phi_E(\zeta)}\right)^{n+1}.$$ Returning to the
Lagrange interpolation polynomial we have
$$\dist_E(f_\zeta,\Cal P_n)\le  ||f_\zeta-L_nf_\zeta||_E =  \max_{z\in E}\left|\frac{\omega_n(z)}{\omega_n(\zeta)
(\zeta-z)}\right|\le  \frac{(n+1)\,c}{\dist\bigl(C(E,\rho),E\bigr)
\dist(\zeta,E)}
\left(\frac{\rho}{\Phi_E(\zeta)}\right)^{n+1}.\qed$$
\enddemo

\proclaim{Lemma 3.4} Assume that a polynomially convex compact set
$E\subset\Bbb C$ admits $\LS(s)$ and $\HCP(k)$ for some $s,k\ge1$,
i.e.\ there exist $a_1,a_2\ge1$ such that for all $z\in E_1$
$$\frac1{a_1}\  \dist(z,E)^s \ \le \
g_E(z) \ \le \ a_2 \, \dist(z,E)^{1/k}.\tag 4$$ Then there exist
$c_0,c_1\ge1$ dependent only on $E$ such that
$$\forall\ \ell\ge1\quad
\forall \ 0<t\le1\quad :\quad \sup_{n\in\Bbb N}\
\frac{n^\ell}{\phi_n(t)} \ \le \ \left(\frac{c_1
\ell}{t^s}\right)^{\ell+c_0} .$$
\endproclaim

\demo{Proof} Fix $\ell\ge1$ and $0<t\le1$.
By Lemma 3.3
for arbitrary $\zeta\in{\text d} E_t$,
$\rho:=\sqrt{\Phi_E(\zeta)}>1$ and $n\in\Bbb N$ we have
$$n^\ell\dist_E(f_\zeta,\Cal P_{n-1}) \! \le \!
\frac{c\,n^{\ell+1}\rho^{-n}}{\dist\bigl(C(E,\rho),E\bigr)  t} \!
\le \! \frac{c}{\dist\bigl(C(E,\rho),E\bigr) t}
\left(\frac{\ell+1}{e \log\rho}\right)^{\ell+1}\!\!\!\!\! \le \!
\frac{c}{\dist\bigl(C(E,\rho),E\bigr)  t}
\left(\frac{2\ell}{g_E(\zeta)}\right)^{\ell+1}\!\!\!,$$ because
for $a,b\!>\!0$ we have \ $\sup_{n>0}n^a  e^{-b n}\!=\!\left(\frac
a{b e}\right)^a\!.$ We combine this with Lemma 3.1 to obtain
$$\frac{n^\ell}{\Phi_n(\zeta)}\le
\bigl( \dist(\zeta,E)+\diam E\bigr) n^{\ell}\dist_E(f_\zeta,\Cal
P_{n-1}) \le \frac{\widetilde c}{\dist\bigl(C(E,\rho),E\bigr)  t}
\left(\frac{2\ell}{g_E(\zeta)}\right)^{\ell+1},$$ where
$\widetilde c:=(1+\diam E) \, c$.  By the above and from
assumption $(4)$,
$$\frac{n^\ell}{\Phi_n(\zeta)}\le
\frac{\widetilde c}t \left(\frac{a_2}{\log \rho}\right)^k
\left(\frac{2\ell}{g_E(\zeta)}\right)^{\ell+1}\frac{\widetilde c}t \left(\frac{2a_2}{g_E(\zeta)}\right)^k
\left(\frac{2\ell}{g_E(\zeta)}\right)^{\ell+1}\!\!\! \le
\frac{\widetilde c}t  \left(\frac{2a_1  a_2}{t^s}\right)^k
\left(\frac{2a_1 \ell}{t^s}\right)^{\ell+1}\!\!\!\le
\left(\frac{c_1 \ell}{t^s}\right)^{\ell+c_0},$$ where $c_0:=k+2$
and $c_1:=2a_1  a_2\widetilde c $ \ depend only on the set $E$.
Finally we conclude that
$$\sup_{n\in\Bbb N}\frac{n^\ell}{\phi_n(t)}\ = \ \sup_{n\in\Bbb N}
\sup_{\zeta\in{\text d} E_t} \frac {n^{\ell}}{\Phi_n(\zeta)} \ \le
\  \left(\frac{c_1 \ell}{t^s}\right)^{\ell+c_0} .\qed$$
\enddemo

\proclaim{Proposition 3.5} For any compact set $E\subset\Bbb C$
and $s\ge1$ the Jackson Property $\JP(s)$ is equivalent to the
following condition:
$$\exists \,\widetilde c,v\ge1\quad \forall\,\ell\ge1\quad
 \forall\,0<t\le1\quad
\forall\, n\in\Bbb N \quad:\quad \frac{n^\ell}{\phi_{n+1}(t)} \ \le
\  \left(\frac{\widetilde c\ell^v}{t^s}\right)^{\ell+\widetilde c}\!\!\!.\tag 5$$
\endproclaim

\demo{Proof} First, observe that we can write equivalently $\ell
\ge 1$ in condition $(1)$ instead of $\ell \in \Bbb N$. Assume
that the set $E$ admits $\JP(s)$ and we shall prove $(5)$. Fix
$0<t\le1$ and arbitrary $\zeta\in{\text d}E_t$. Obviously,
$f_\zeta\in\Cal H^{\infty}(E_\delta)$ for $\delta:=\frac t2$.
$\JP(s)$ implies that
$$n^\ell  \dist_E(f_{\zeta},\Cal P_n)\le
|f_{\zeta|E}|_\ell\le
\left(\frac{c\ell^v}{\delta^s}\right)^{\ell+c}
||f_\zeta||_{E_\delta} = \left(\frac{2^s
c\ell^v}{t^s}\right)^{\ell+c} \cdot \frac2t\le \left(\frac{2^{s }
c\ell^v}{t^s}\right)^{\ell+c+1}$$ for all $\ell\ge1$ and
$n\in\Bbb N$. By Lemma 3.1, we obtain
$$\frac{n^\ell}{\Phi_{n+1}(\zeta)}\ \le \
\bigl(\dist(\zeta,E)+\diam E\bigr) \, n^\ell \dist_E(f_\zeta,\Cal
P_n) \ \le \ \left(\frac{\widetilde
c\ell^v}{t^s}\right)^{\ell+\widetilde c},
$$
where $\widetilde c:=\max\{(1+\diam E)\,2^{s }c,c+1\}$.
Therefore, since $\zeta\in{\text d}E_t$ was
arbitrary, we conclude that
$$\frac{n^\ell}{\phi_{n+1}(t)} \ = \ \sup_{\zeta\in{\text d}E_t}\frac{n^\ell}{\Phi_{n+1}(\zeta)} \
\le\  \left(\frac{\widetilde c\ell^v}{t^s}\right)^{\ell+\widetilde c}$$ and $(5)$ is proved.

In order to show $\JP(s)$ assuming $(5)$, fix $\ell\ge1$,
$0<\delta\le1$ and $f\in\Cal H^\infty(E_\delta)$. We apply
Prop.~3.2 with $b:=\frac12$ and the assumption with  $t:= \delta/2$
to obtain for any $n\in\Bbb N$
$$ n^\ell  \dist_E(f,\Cal P_n)\ \le \ \frac{2c n^\ell
\|f\|_{E_\delta}}{\delta^2 \phi_{n+1}(\delta/2)} \ \le \
\frac{2c}{\delta^2} \left(\frac{2^s  \widetilde c\ell^v}
{\delta^s}\right)^{\ell+\widetilde c} \|f\|_{E_\delta} \
\le \ \left(\frac{2^{s }c \, \widetilde c\ell^v}
{\delta^s}\right)^{\ell+\widetilde c+2} \|f\|_{E_\delta}
.$$ From this it follows that for $c_0:=\max\{1+2^{s } c \,\widetilde c,\widetilde c+2\}$ we have
$$|f_{|E}|_\ell \ = \ \|f\|_E+\sup_{n\in \Bbb N}n^\ell  \dist_E(f,\Cal P_n)\ \le
\ \left(\frac{c_0\ell^v}{\delta^s}\right)^{\ell+c_0}
||f||_{E_\delta}.\qed$$   \enddemo

\proclaim{Theorem 3.6} Let $s\ge1$. Any polynomially convex
compact set $E\subset\Bbb C$ admitting $\LS(s)$ and $\HCP$, admits
$\JP(s)$ with $v=1$.
\endproclaim

\demo{Proof} This is an immediate consequence of Lemma 3.4 and
Prop.~3.5. \qed
\enddemo

The closest we could find in the literature was an estimate
equivalent to $\JP(1)$ with $v=1$ and $c\ge2$, proved for all simply
connected bounded regions with boundaries that are Jordan curves
of class $\Cal C^{1+\varDelta}$ \cite{17, lemma 4}.

Note that as a simple corollary of Theorem 3.6, we can obtain
$\JP(1)$ for a disk $E=B(0,r)$, because in this case we have
$\Phi_E(z)=|z|/r$.

\proclaim{Proposition 3.7} For any compact set $E\subset\Bbb C$
and $s'>s\ge1$  we have
$$\JP(s)\implies \LS(s').$$
\endproclaim

\demo{Proof}
% Assume that $E$ admits $\JP(s)$.
By Prop.~3.5,
for arbitrary $t\!\in\!(0,1]$, $\zeta\!\in\!{\text
d}E_t$, $n\!\in\!\Bbb N$ and $\ell\!\ge\!1$
we get
$$g_E(\zeta) \ = \ \log\Phi_E(\zeta) \ \ge \
\log\root {n+1}\of{\Phi_{n+1}(\zeta)} \ \ge \ \log\root
{n+1}\of{\phi_{n+1}(t)} \ \ge \ \frac1{n+1} \ \log\left(n^\ell
\left(\frac{t^s}{\widetilde c\ell^v}\right)^{\ell+\widetilde c}\right).$$
Specifically, by
taking $n\!\in\!\Bbb N \cap \left[e \!\left(\frac{\widetilde
c\ell^v}{t^s}\right)^{\!1+\widetilde c/\ell}\!\!\!,
1\!+\!e\! \left(\frac{\widetilde
c\ell^v}{t^s}\right)^{\!1+\widetilde c/\ell}\right)$ we obtain
$$g_E(\zeta) \ \ge \ \frac{\log e^\ell}{2+ e\!
\left(\frac{\widetilde c\ell^v}{t^s}\right)^{1+\widetilde
c/\ell} } \ \ge
\ \frac{\ell}{2+e \, (\widetilde
c\ell^v)^{1+\widetilde c/\ell}}\ t^{s (1+\widetilde c/\ell)}
\ = \ \frac{\ell}{2+e\, (\widetilde c\ell^v)^{1+\widetilde
c/\ell}} \ \dist(\zeta,E)^{s (1+\widetilde c/\ell)}.$$
If we take $\ell$ sufficiently large then we obtain $\LS(s')$ for any $s'>s$.
% If we take the limit $\ell\rightarrow\infty$ in the last estimate, then we obtain
% $$\gather
% g_E(\zeta) \ \ge \ \lim_{\ell\rightarrow\infty}
% \frac{\ell\,\dist(\zeta,E)^{s (1+\widetilde c_0/\ell)}}{2+e\,
% (\widetilde c_1\, \ell)^{\,1+\widetilde c_0/\ell}} \
%  = \ \frac{1}{e \widetilde
% c_1} \ \dist(\zeta,E)^s.\qed\endgather$$
\qed
\enddemo
Note that if we have $\JP(s)$ with $v=1$ in the assumption of the last proposition,
then we can conclude $\LS(s)$ rather than $\LS(s')$ for any $s'>s$,
by simply taking the limit for $\ell\rightarrow+\infty$ in the last inequality of the proof.

\vskip 5mm

\head 4. Remarks and additional results.
\endhead

\proclaim{Proposition 4.1}
Every compact set $E\subset\Bbb R$
% , which admits $\HCP$, also
admits $\JP(1)$.
\endproclaim

To the extent that the set $E$ admits $\HCP$,
this proposition is a simple corollary
of the main theorem and the fact that
every compact set in $\Bbb R$ admits $\LS(1)$.
We leave the proof of the general case to the reader.
Hint: apply the classical Jackson inequality and
appropriate cut-off functions to estimate the quotient norms.
Our best estimate gives $v=6$.

\vskip 2mm

One may ask whether the Jackson Property is maintained after combining two sets into one,
or separating one set into two distinct subsets.
Before answering this question we need to do some preparations.

\proclaim{Lemma 4.2} Assume that the compact set $E\subset\Bbb C$
is the sum of two polynomially convex, disjoint compact subsets,
i.e.\ $E=A\cup B$, $A=\hat A$, $B=\hat B$, $A\cap B=\emptyset$.
Assume also that the subset $A$ is non-polar, i.e.\ $\capa A>0$.
Then for any function $f\in\Cal C(E)$ such that $f_{|A}\in s(A)$
and $f_{|B}\equiv 0$, we have $f\in s(E)$ and furthermore we can
estimate its \Zerner norms on the set $E$ by its \Zerner norms on
the subset $A$ as follows:
$$\forall\ \ell\ge1\quad :\quad
|f|_\ell \ \le \ (c \ell)^\ell \, |f_{|A}|_\ell,$$ where the
constant $c\ge1$ depends only on the subsets $A$ and $B$. Note
that these are two different \Zerner norms and only the domain of
the function indicates which norm is meant.
\endproclaim

\demo{Proof} Like in the proof of \cite{18, Th.$\,$1} we consider
$$\chi_B (z)=\cases
0&\quad\text{for \ }z\in A,\\
1&\quad\text{for \ }z\in B,\endcases$$ and we note that this
characteristic function can be extended holomorphically so that
$\chi_B\!\!\in\!\Cal
H^\infty\!\left(\!\widehat{C(E,2\rho)\!}\right)$ for some
$\rho>1$. By \cite{13, chap.II \S 3A Th.1}, there exists a
constant $M\ge1$ such that
$$\forall \, n\in\Bbb N\quad: \quad
\dist_E(\chi_B,\Cal P_n) \ \le \ \frac M{\rho^n}.$$ Put
$x:=||\Phi_A||_B$ and note that $1<x<+\infty$, because the subset
$A$ is non-polar and both subsets $A$ and $B$ are compact.
Therefore, we can find an integer $k\in\Bbb N$ such that
$t:=\frac{\rho^k}x>1$.

Now fix an arbitrary function $f\in\Cal C(E)$, such that
$f_{|A}\in s(A)$ and $f_{|B}\equiv 0$, and also fix a number
$\ell\ge1$. Find two sequences of polynomials of best
approximation for the functions $f_{|A}$ and $\chi_B$ on the sets
$A$ and $E$ respectively, i.e.\ $p_n,q_n\in \Cal P_n$,
$\|f-p_n\|_A=\dist_A(f,\Cal P_n)$ and
$\|\chi_B-q_n\|_E=\dist_E(\chi_B,\Cal P_n)$ for each $n\in \Bbb
N_0$. Clearly, $\|p_n\|_A \le \|f\|_A+\|f-p_n\|_A \le 2\|f\|_A$.
By the definition of Siciak's extremal function, we see that
$$\|p_n\|_B \ \le \
\|\Phi_A\|_B^n \ \|p_n\|_A \ = \ x^n \|p_n\|_A \ \le  \ 2x^n
\|f\|_A.$$

For each $n\in \Bbb N_0$ we put
$$r_n(z) \ := \ p_n(z) \,\bigl(1-q_{k  n}(z)\bigr)
$$
so that $r_n\in\Cal P_{(k+1)  n}$. This way we obtain
$$\|f-r_n\|_A \ \le \ \dist_A(f,\Cal P_n)+\|p_n\|_A \dist_A(\chi_B,\Cal
P_{k n}) \ \le \ \frac{|f_{|A}|_\ell}{n^\ell}+2 \|f\|_A \frac
M{\rho^{k n}} \ < \ \frac{|f_{|A}|_\ell}{n^\ell}+\frac {2M}{t^n}
\|f\|_A,$$ $$\|f-r_n\|_B \ = \ \|r_n\|_B \ \le \ \|p_n\|_B\
\|\chi_B-q_{k n}\|_B \ \le \ \|p_n\|_B \dist_E(\chi_B,\Cal P_{k
n})\ \le \ 2x^n \|f\|_A \frac M{\rho^{k  n}} \ =  \ \frac
{2M}{t^n} \|f\|_A.$$ This then leads us to
$$n^\ell \dist_E(f,\Cal P_{(k+1)  n})\le n^\ell \|f-r_n\|_E\le
|f_{|A}|_\ell+2M  \frac {n^\ell}{t^n} \|f\|_A \le|f_{|A}|_\ell+2M
\left(\frac {\ell}{e \log t}\right)^\ell \|f\|_A\le (\widetilde c
\ell)^\ell |f_{|A}|_\ell,$$ where the constant $\widetilde
c:=1+\frac{2M}{e \log t}$ depends on the subsets $A$ and $B$ but
not on the choice of the function $f$ and the number $\ell$.

Finally, for arbitrary $n\!\in\!\Bbb N_0$ we can find
$N\!\in\!\Bbb N_0$ such that $(k+1)  N\!\le\! n\!<\!(k+1) (N+1)$
to conclude~that
$$n^\ell \dist_E(f,\Cal P_n)\le \bigl((k+1) (N+1)\bigr)^\ell
\dist_E(f,\Cal P_{(k+1)  N}) \le(4k)^\ell  N^\ell  \dist_E(f,\Cal
P_{(k+1)  N})\le (4k \widetilde c \ell)^\ell |f_{|A}|_\ell, $$
$$|f|_\ell = \|f\|_E+\sup_{n\in\Bbb N}n^\ell \dist_E(f,\Cal
P_n)\le \|f\|_A+ (4k \widetilde c \ell)^\ell |f_{|A}|_\ell\le (c
\ell)^\ell |f_{|A}|_\ell< +\infty,$$ where the constant $c:=1+4k
\widetilde c$ also depends only on the subsets $A$ and $B$. \qed
\enddemo

\proclaim{Corollary 4.3} Assume that the compact set $E\subset\Bbb
C$ is the sum of two polynomially convex, disjoint, non-polar
compact subsets, i.e.\ $E=A\cup B$, $A=\hat A$, $B=\hat B$, $A\cap
B=\emptyset$, $\capa A>0$ and $\capa B>0$. Then for any function
$f\in\Cal C(E)$ such that $f_{|A}\in s(A)$ and $f_{|B}\in s(B)$,
we have $f\in s(E)$ and furthermore we can estimate its \Zerner
norms on the set $E$ by its \Zerner norms on the subsets $A$ and
$B$ as follows:
$$\forall\,\ell\ge1\quad :\quad
|f|_\ell \ \le \ (c \ell)^\ell
\left(|f_{|A}|_\ell+|f_{|B}|_\ell\right),$$ where the constant
$c\ge1$ depends only on the subsets $A$ and $B$. Note that these
are three different \Zerner norms and only the domain of the
function indicates which norm is meant.
\endproclaim

\demo{Proof} We put $\widetilde f^A:=\chi_A \cdot f$ and
$\widetilde f^B:=\chi_B\cdot f$ \ so that $\widetilde
f^A,\widetilde f^B\in\Cal C(E)$. We apply Lemma 4.2 to obtain
$$|f|_\ell \!\le
|\widetilde f^A|_\ell+|\widetilde f^B|_\ell \!\le (c_A \ell)^\ell
|\widetilde f^A_{|A}|_\ell + (c_B \ell)^\ell |\widetilde
f^B_{|B}|_\ell \! = (c_A \ell)^\ell |f_{|A}|_\ell+ (c_B \ell)^\ell
|f_{|B}|_\ell\!\le (c \ell)^\ell
\left(|f_{|A}|_\ell\!+|f_{|B}|_\ell\right)\!<\!+\infty$$ for any
$\ell\ge1$  with the constant $c:=\max\{c_A,c_B\}$ depending only
on the sets $A$ and $B$. \qed
\enddemo

\proclaim{Proposition 4.4} Assume that the compact set
$E\subset\Bbb C$ is the sum of two polynomially convex, disjoint
compact subsets, i.e.\ $E=A\cup B$, $A=\hat A$, $B=\hat B$ and
$A\cap B=\emptyset$. If the set $E$ admits $\JP(s)$ with some
$s\ge1$, then both subsets $A$ and $B$ admit $\JP(s)$.

Conversely, if both subsets $A$ and $B$ are additionally non-polar
and they both admit $\JP(s)$ with some $s\ge1$, then the set $E$
admits $\JP(s)$.
\endproclaim

\demo{Proof}
In order to prove the first assertion,
we note that the Jackson Property is invariant to an affine change of variable
and therefore if necessary we can blow these sets up so that $\dist(A,B)>2$.
This way the intersection of the neighborhoods $A_1$ and $B_1$ of the sets $A$ and $B$,
respectively,
is empty.
Next we apply Prop.~3.5 to get condition $(5)$ for the set $E$.
The extremal functions $\Phi_n$
of the subsets $A$ and $B$ are bounded below by the respective
extremal functions of the set $E$ and this way we obtain the
condition $(5)$ for the sets $A$ and $B$.
Finally, we apply
Prop.~3.5 again to conclude that they too admit $\JP$ with the same
coefficients.

The second assertion follows straight from Corollary 4.3 and the
definition of the Jackson Property. Indeed, for arbitrary
$\ell\ge1$, $0<\delta\le1$ and $f\in\Cal H^\infty(E_\delta)$ we
have $f\in\Cal H^\infty(A_\delta)$, $f\in\Cal H^\infty(B_\delta)$
and
$$\gather
|f_{|E}|_\ell\le (c \ell)^\ell
\left(|f_{|A}|_\ell+|f_{|B}|_\ell\right)\le\\
\le(c \ell)^\ell \left(
\left(\frac{c_A\ell^v}{\delta^s}\right)^{\ell+c_A}
||f||_{A_\delta}+
\left(\frac{c_B\ell^v}{\delta^s}\right)^{\ell+c_B}
||f||_{B_\delta}\right)\le
\left(\frac{\widetilde c\ell^{v+1}}{\delta^s}\right)^{\ell+\widetilde c}
||f||_{E_\delta},
\endgather$$
where $\widetilde c:=2c \max\{c_A,c_B\}$.
% and $c_0:=\max\left\{c_0^A,c_0^B\right\}$.
\qed
\enddemo

%\vskip 3mm

\noindent {\bf Remark 4.5.}
We close this paper by offering three
open problems for further research:

\vskip 2mm

$\bullet$ \
The proof of Lemma 3.4 applies the assumption of $\HCP$ only in order to make sure
that the level sets of the extremal function do not come too close to the compact set $E$.
The coefficient in $\HCP(k)$ has no meaningful impact on the coefficients of the Jackson Property,
suggesting that we may have used a sledge-hammer to crack a nut.
% However multiple attempts at weakening this assumption have all failed.
Specifically,
due to the intended application of the Jackson Property,
it would be interesting to know whether it is sufficient to assume $\GMI$
instead of $\HCP$ (which implies $\GMI$)?
It should be noted though that Lemma 3.3 assumes L-regularity,
which is guaranteed by $\HCP$,
but it is still not known whether all compact subsets of the complex plane admitting $\GMI$ are L-regular.
In the real case this follows from the combination of \cite{7} and \cite{4}.

\vskip 1mm

$\bullet$ \
The characterization of compact sets $E\subset\Bbb C$,
for which $\Cal A^\infty(E)_{|E}=s(E)$,
also remains an open problem,
especially for totally disconnected sets.
Siciak proved this property for simply connected H\"older domains,
i.e.\ admitting $\LS$ \cite{22, Th.1.10}.
More recently,
Gendre proved the same for every compact set
$E\subset\Bbb C^N$ that is Whitney $1$-regular and admits $\HCP$ as well as $\LS$
\cite{14, Cor.\ 7}.

\vskip 1mm

$\bullet$ \
Finally we had a good look at the Wiener type characterization given by Carleson and Totik
for pointwise H\"older continuity of Green's functions.
Their Wiener type criterion (i.e.\ lower bounds for capacities) introduced in \cite{8} implies $\HCP$,
but in order to assert the converse they needed an additional assumption,
i.e.\ either a (geometric) cone condition or a quantitative (capacity) condition (upper bounds for capacities).
The examples given in the Introduction above suggest that both those conditions could be special cases of $\LS$.
It is worth investigating whether $\HCP$ in conjunction with $\LS$ is sufficient to assert
the Wiener type criterion proposed by Carleson and Totik.

\enddocument